\documentclass[a4paper]{article}
\usepackage{amsmath,amsthm,amstext,amsbsy,amssymb,amsfonts,amscd}
\usepackage{aeguill}

\newtheorem{Th}{Theorem}[]
\newtheorem{Lem}[Th]{Lemma}
\newtheorem{Prop}[Th]{Proposition}
\newtheorem{Cor}[Th]{Corollary}

\newtheorem{Def} [Th]{Definition}

\newtheorem{DProp}[Th]{Definition \& Proposition}

\newtheorem{Rem}[Th]{Remark}
\newtheorem{alg}[Th]{Algorithm}
\def\vpX{\vrule height 12pt depth 4pt width 0pt}
\newcommand{\Algo}[5]
            {
            \begin{alg}[#1] \label{#2}{$\;$}\rm
                \\
\mbox{\enspace}
                \rlap{\rm Input:\enspace}\phantom{\rm Output:\enspace}
\vtop{\hsize=130mm\noindent#3}
                \\
\mbox{\enspace}
                {\rm Output:\enspace}
\vtop{\hsize=130mm\noindent#4}
\parskip0pt
\begin{list}{}{\setlength{\leftmargin}{0pt}}
\item                #5 
\end{list}
\parskip6pt
            \end{alg}
            \goodbreak}
\def\itemb{\item[$\bullet$]}

\def\goth{\mathfrak}
            
\def\F2{\mathbb{F}_2}   
\def\Fp{\mathbb{F}_p}        
\def\F{\mathbb{F}}   
\def\CC{\mathbb{C}}   
\def\Q{\mathbb Q}
\def\Qp{\mathbb{Q}_p}        
\def\Zl{\mathbb{Z}_\ell}        
\def\Zp{\mathbb{Z}_p}        
\def\Z2{\mathbb{Z}_2}
\def\Z{\mathbb Z}

\def\U{\mathcal  U}
\def\J{\mathcal  J}     
     
\def\R{\mathcal  R}
\def\Dl{\mathcal  D\ell}        
\def\Pl{\mathcal  P\ell}        
\def\Cl{\mathcal  C\ell}
\def\Id{\mathcal Id}    
\def\Pr{\mathcal  Pr}
\def\E{\mathcal  E}

\def\q{{\goth q}}               
\def\s{{\goth s}}               
\def\p{{\goth p}}               
\def\l{{\ell}}
\def\d{{\goth d}}               
\def\a{{\goth a}}
\def\b{{\goth b}}
\def\g{{\goth g}}

\def\wi{\widetilde}             
\def\ov{\overline}

\def\div{\wi {\operatorname{div}}}      
\def\deg{\operatorname{deg}}
\def\rank{\operatorname{rank}}
\def\Gal{\operatorname{Gal}}    
\def\Log{\operatorname{Log}}
\def\Ker{\operatorname{Ker}}    
\def\Coker{\operatorname{Coker}}
\def\norm{\operatorname{N}}

\def\Fp{\mathbb{F}_p}
\def\OOk{{\mathcal O}_K}

\def\OOp{{\mathcal O}_{\p}} 

\def\Fe{{\mathcal F}_e}
\def\enul{{e_0}}
\def\munul{{\mu_0}}

\begin{document}

\title{\bf A new Algorithm for the Computation of logarithmic
$\ell$-Class Groups of Number Fields\footnote{Experimental. Math. {\bf 14} (2005), 67--76.}} 

\author{{\normalsize Francisco {\sc Diaz y Diaz}, Jean-Fran\c{c}ois 
{\sc Jaulent}, Sebastian {\sc Pauli},}  
\and\normalsize {Michael {\sc Pohst}, Florence {\sc Soriano-Gafiuk}}}

\date{}
\maketitle

\bigskip

{\small
\noindent {\bf Abstract.} We present an algorithm for the computation of logarithmic $\l$-class groups of number fields.  Our principal motivation is the effective determination of the $\l$-rank of the wild kernel in the $K$-theory of number fields.\medskip

\noindent {\bf R\'esum\'e}. Nous d\'eveloppons un nouvel algorithme pour calculer le $\l$-groupe $\widetilde \Cl_F$ des classes logarithmiques d'un corps de nombres $F$. Il permet en particulier de d\'eterminer effectivement le $\l$-rang du noyau sauvage $W\!K_2(F)$ du groupe $K_2(F)$.
}

%%%%%%%%%%%%%%%%%%%%%%%%%%%%%%%%%%%%%%%%%%%
%%%%%%%%%%%%%%%%%%%%%%%%%%%%%%%%%%%%%%%%%%%

\section{Introduction}

A new invariant of number fields, called group of logarithmic classes,
was introduced by J.-F. Jaulent in 1994 \cite{J1}.
The interest in the  arithmetic of logarithmic classes is because
of its applicability in $K$-Theory. 
Indeed this new group of classes
is revealed to be mysteriously related to the wild kernel in the  
$K$-Theory for number fields.
The new approach to the wild kernel is so attractive since
the arithmetic of logarithmic classes is very efficient.
Thus it provides an algorithmic and original study of the wild kernel.
A first algorithm for the computation of the group of logarithmic 
classes of a number field $F$ was developed 
by F. Diaz y Diaz and F. Soriano in 1999 \cite{DS}.
We present a new much better performing algorithm, which also eliminates the 
restriction to Galois extensions.

Let $\l$ be a prime number.
If a number field $F$ contains the $2\l$-th roots of unity then the wild kernel
of $F$ and its logarithmic $\l$-class group have the same $\l$-rank.
If $F$ does not contain the $2\l$-th roots of unity the arithmetic 
of the logarithmic classes still yields the $\l$-rank of the the wild kernel.
More precisely:
\begin{itemize}
\item
If $\l$ is odd \cite{JS,So} we consider $F':=F(\zeta_\l)$, where $\zeta_\l$ is the
$\l$-th root of unity, and use classic techniques from the theory of semi-simple algebras.
\item
If $\l=2$ \cite{JS2} we introduce a new group, which we call the $\l$-group of 
of the positive divisor classes and which can be constructed from the $\l$-group of
logarithmic classes.
\end{itemize}

In the present article we consider the general situation where $F$ is a number field which does not necessarily contain the $2\l$-th roots of unity.

%%%%%%%%%%%%%%%%%%%%%%%%%%%%%%%%%%%%%%%%%%%
%%%%%%%%%%%%%%%%%%%%%%%%%%%%%%%%%%%%%%%%%%%

%%%%%%%%%%%%%%%%%%%%%%%%%%%%%%%%%%%%%%%%%%%
%%%%%%%%%%%%%%%%%%%%%%%%%%%%%%%%%%%%%%%%%%%

\section{The theoretical  background}

%%%%%%%%%%%%%%%%%%%%%%%%%%%%%%%%%%%%%%%%%%%
%%%%%%%%%%%%%%%%%%%%%%%%%%%%%%%%%%%%%%%%%%%

This section is devoted to the introduction of the main notions of logarithmic 
arithmetic. We also review  the facts that are of interest for our purpose. 
We do not attempt to give a fully detailed account of the logarithmic 
language. Most proofs may  be found in \cite{J1}, pp 303-313.

%%%%%%%%%%%%%%%%%%%%%%%%%%%%%%%%%%%%%%%%%%%

\subsection{Review of the main logarithmic objects}

%%%%%%%%%%%%%%%%%%%%%%%%%%%%%%%%%%%%%%%%%%%

For any number field  $F$, let  $\J_F$ be the 
{\it $\l$-adified  group of id\`eles} of $F$, i.e. the restricted product

\[
{\J_F\, =\, \prod_{\p}^{res}\R _{\p}}
\]

\noindent of the  $\l$-adic  compactifications  $\R _{\p}\, =\ \varprojlim
\ {F_{\p} ^\times }/{F_{\p}^\times }^{\l^n}$
of the   multiplicative groups of the  completions of  $F$ at each $\p$.
For each  finite  place  $\p$ the subgroup ${\wi\U }_\p$
of $\R _\p$ of the  cyclotomic norms  (that is to say the elements 
of  $\R _\p$ which  are  norms at any finite step of the 
local cyclotomic $\Z_\l$-extension  $F_{\p}^c/F_{\p}$)
will be called {\it the group of  logarithmic units} of $F_{\p}$.  
The product
\[
{\wi\U_{F}\, =\, \prod_\p \ \wi\U _{\p}}
\]
is denominated   the  {\it group of idelic logarithmic units}; 
it happens to be the kernel of the {\it logarithmic valuations}
\[
{\wi v_\p\ |  \ x\, \mapsto\, -\frac{\Log_\ell\, (N_{F_\p/\Q_p}
(x))}{\deg_F \p} ,}
\]
defined on the   $\R _{\p}$ and $\Z _{\l}$-valued. These are 
 obtained by taking the Iwasawa logarithm of the norm of 
 $x$ in the local extension  $F_{\p}/\Q _p$ with a normalization factor $
\deg _F \, \p$ whose precise definition is given in the next subsection.

The quotient $\Dl _F\, =\, \J _F /\wi{\U}_F$
is the  $\l$-group of {\it  logarithmic divisors} of  $F$;  via the 
logarithmic valuations  $\wi v_{\p}$, it 
may be identified with the free  $\Z _{\l}$-module generated by the prime 
ideals of $F$
\[
{\Dl _F\, =\, \J _F /\wi{\U}_F\, =\, \oplus _{\p}\ \Zl \, \p.}
\]

The  {\it degree} of a  logarithmic divisor 
$\d  =\sum _{\p} n_{\p} \, \p$ is then defined by
\[
\deg _F\, (\sum _{\p} n_{\p} \, \p )\, =\, \sum _{\p} 
n_{\p}\deg _F\, \p,
\]
\noindent inducing a $\Z_\l$-valued $\Z_\l$-linear map on the class group of 
logarithmic divisors.
The  logarithmic divisors of   degree zero form a subgroup of ${\Dl} _F$
denoted by \smallskip

\centerline{$\wi{\Dl} _F\, =\, \{\, {\d} \in {\Dl }_{F}\, |  \, \deg _F\,{\d}  \, =\, 0\, \}$.}
\smallskip

The image of  the map   $\div _F$ defined via the set of 
logarithmic valuations from the principal id\`ele  subgroup\smallskip

\centerline{$\R _F\, =\, \Zl \otimes _\Z F^\times$}\smallskip

\noindent of  $\J _F$ to $\wi\Dl _F$ is a subgroup denoted by $\wi \Pl  _F$, 
which will be referred to as the subgroup  of {\it  principal logarithmic divisors}.
The  quotient 
$$
\widetilde \Cl _F\, =\, \widetilde{\Dl} _F / \widetilde {\Pl } _F
$$
is, by  definition, the  {\it $\l$-group of logarithmic classes} of $F$.
And the kernel\smallskip

\centerline{$\wi \E _F\, =\, \R _F \cap \wi \U_F$}\smallskip

\noindent of the morphism $\div _F$ from $\R_F$ in  
${\wi \Dl} _F$ is the group of global {\it logarithmic units}.

%%%%%%%%%%%%%%%%%%%%%%%%%%%%%%%%%%%%%%%%%%%

\subsection{Logarithmic ramification and $\l$-adic-degrees.}\label{subseclogram}

%%%%%%%%%%%%%%%%%%%%%%%%%%%%%%%%%%%%%%%%%%%

Next we review the basic notions of  the logarithmic 
ramification, which  mimic, as  a rule, the classical ones. 

Let $L/F$ be any
finite extension of number fields.  Let $p$ be a prime number. 
Denote by $\widehat{\mathbb Q}^c_p$ 
the cyclotomic $\widehat{\mathbb Z}$-extension of ${\mathbb Q}_p$,
that is to say the compositum of all cyclotomic $\mathbb Z_q$-extension of  
$\mathbb Q_p$ on all prime numbers $q$.
Let $\goth p$ be a prime of $F$ above $(p)$ and ${\goth P}$
a prime of $L$ above $\goth p$. The logarithmic ramification (resp. inertia)
 index  $\wi{e}(L_{\goth P}/F_\p)$ (resp. $\wi{f}(L_{\goth P}/F_\p)$)
is defined to be the relative degree 
$$
\wi{e}(L_{\goth P}/F_\p) = [L_{\goth P} : L_{\goth P}\cap \widehat
{\mathbb Q}^c_p F_{\goth p}] \quad ({\rm resp. \quad } \wi{f}(L_{\goth P}/F_\p
)= [ L_{\goth P} \cap \widehat{\mathbb Q}^c_p F_{\goth p}:F_{\goth p}]). 
$$
As a  consequence, $L/F$ is logarithmically unramified  at 
${\goth P}$, i.e. $\widetilde e(L_{\goth P}/F_{\goth p})=1$, 
if and only if  $L_{\goth P}$ is contained in  the 
cyclotomic extension of $F_{\goth p}$.  Moreover, for any  
$q\ne p$  the classical and the logarithmic 
indexes  have the same  $q$-part (see theorem \ref{dprop5}). 
Hence  they are equal as soon as 
 $p\nmid [F_{\goth p} : \mathbb Q_p]$.

As usual, in the special case $L/F=K/\Q$, the absolute logarithmic indexes 
of a finite place $\p$ of $K$ over the prime $p$ are just denoted by $\wi e_\p$
and $\wi f_\p$. With these notations, the $\ell$-adic degree of $\p$  
is defined by the formula:
\[
\deg_K \p = \wi f_\p \deg_\l p \quad {\rm with }\quad \deg_\l p =\left\{\aligned
&\Log_\ell p \  &{\rm for}\ p \ne \ell ; \\
&\ell \   &{\rm for}\ p =\ell \ne 2; \\
&4 \    &{\rm for}\ p =\ell = 2.
\endaligned
\right.
\]

The extension and norm maps  between groups
of divisors, denoted by 
${\iota}_{L/F}$ and ${N}_{L/F}$ respectively, have their logarithmic 
counterparts,  $\wi \iota_{L/F}$ and $\wi N_{L/F}$ respectively. To be more 
explicit,   $\widetilde{\iota}_{L/F}$ is defined on every finite place 
$\p$ of $F$ 
by\smallskip

\centerline{$\wi {\iota}_{L/F}({\goth p})=\sum_{{\goth P|\goth p}}
\wi e_{L_{\goth P}/F_{\goth p}}{\goth P}$ ,}\smallskip

\noindent while  $\wi {N}_{L/F}$ is defined on  all  $\goth P$ lying above 
$\goth p$ by\smallskip
  
\centerline{$\wi{ N}_{L/F}({\goth P})={\wi f}_{L_{\goth P}/F_\p} \ \p$ .}\smallskip

\noindent These applications are  compatible with the usual extension and 
norm maps defined between  $\mathcal R_L$ and  $\mathcal R_F$, in the 
sense that they sit inside the commutative diagrams:
$$
\CD 
\mathcal R_L     @>\div_L>> \wi{\Dl}_{L}    @>\deg_L>>      \Zl @. \qquad
\qquad  @. \R_L     @>\div_L>>  \wi{\Dl}_{L}    @>\deg_L>>      \Zl\\ 
@VV{N}_{L/F}V        @VV\wi { N}_{L/F}V          @\vert\text{\rm and}@.  
@AA\widetilde\iota_{L/F}A          @AA\widetilde\iota_{L/F}A       
@AA[L:F]A \\    
\mathcal R_F     @>\div_{F}>> \widetilde{\mathcal D\l}_{F}   @>\deg_F>>    
\mathbb Z_\ell @. \qquad \qquad @. \mathcal R_F     @>\div_{F}>>    
\widetilde {\mathcal D\l}_{F}    @>\deg_F>>     \mathbb Z_{\ell}.\\
\endCD
$$
When $L/F$ is a Galois extension with Galois group $\text{Gal}(L/F)$, one 
deduces from the very definitions the unsurprising and obvious properties:
\smallskip

\centerline{$\widetilde{N}_{L/F}\circ \wi \iota_{L/F} =[L:F]$ \quad and \quad $
\wi \iota_{L/F}\circ \wi N_{L/F}=\sum_{\sigma\in\Gal(L/F)} \sigma$ .}
\bigskip

%%%%%%%%%%%%%%%%%%%%%%%%%%%%%%%%%%%%%%%%%%%

\subsection{Ideal theoretic description of logarithmic classes.}\label{subsecideal}

%%%%%%%%%%%%%%%%%%%%%%%%%%%%%%%%%%%%%%%%%%%

By the weak density theorem every class in $\J_F/\wi\U_F\R_F$ may 
be represented by an idele with trivial components at the $\ell$-adic
places, that is to say that every class in $\Dl_F/\Pl_F$ comes from a 
$\ell$-divisor $\d = \sum_{\p\nmid\ell}\ \alpha_\p \ \p$.

The canonical map from $\R_F$ to $\Dl_F$ maps
$a\in\R_F$ to $\div_F(a)=\sum_\p \wi v_\p(a)\p$.
Now for each finite place $\p\nmid\ell$, the quotient $\wi e_\p / e_\p =
f_\p / \wi f_\p$ of the classical and logarithmical indexes associated 
to $\p$ is a $\ell$-adic unit (theorem \ref{dprop5}), 
say $\lambda_\p$ (which is 1 for almost
all $\p$), and one has the identity:\smallskip

\centerline{$\wi v_\p = \lambda_\p v_\p$}\smallskip

\noindent between the logarithmic and the classical valuations. So
every $\ell$-divisor $\d$ comes from a $\ell$-ideal $\a$ by the 
formula:\smallskip

\centerline{$\a = \prod_{\p\nmid\ell}\p^{\alpha_\p} \mapsto \d_F(\a) = 
\sum_{\p\nmid\ell} \lambda_\p\ \alpha_\p \ \p$ .}\smallskip

This gives the following ideal theoretic description of logarithmic 
classes:

%Definition & Proposition 1
\begin{DProp}\label{defpropcl} 
Let $\Id_F = \{\a = \prod_{\p\nmid\ell}\p^{\alpha_\p}  \}$ 
be the group of $\ell$-ideals, 
$\wi\Id_F = \{\a \in \Id_F| \deg_F\d_F(\a)=0 \}$ 
be the subgroup of $\ell$-ideals of degree 0
and $\wi\Pr_F = \{\prod_{\p\nmid\ell}\p^{v_\p(a)} | \wi v_\p(a)=0 \ \forall 
\p \mid \ell\}$ the subgroup of principal
$\ell$-ideals generated by principal ideles $a$ having logarithmic 
valuations 0 at every $\ell$-adic places. Then one has:\smallskip

\centerline{$\wi\Cl_F \simeq \wi\Id_F / \wi\Pr_F$}
\end{DProp}

%\noindent{\it Proof.} 
\begin{proof}
As explained above, the surjectivity follows from
the weak approximation theorem. So let us consider the kernel of 
the canonical map $\phi_F: \ \wi\Id_F \mapsto \wi\Cl_F$. Clearly we have:
$\ker \phi_F = \{\a \in \wi \Id_F \ | \ \exists a \in \R_F \quad \d_F(\a)= \div_F(a)\}$. The
condition $\d_F(\a)= \div_F(a)$ with  $\a \in \wi \Id_F$ implies 
$\wi v_\p(a)=0 \ \forall \p \mid \ell$; 
and then gives 
$(a) \in \wi \Pr_F$ 
as expected.
\end{proof}
%\medskip

The generalized Gross conjectures (for the field $F$ and the prime
$\ell$) asserts that the logarithmic class group $\wi\Cl_F$ is finite
(cf. \cite{J1}). This conjecture, which is a consequence of the $p$-adic
Schanuel conjecture was only proved in the abelian case and a few 
others (cf. \cite{FG,J5}). Nevertheless, since $\wi\Cl_F$ is a $\Zl$-module 
of finite type (by the $\ell$-adic class field theory), the Gross conjecture
just claims the existence of an integer $m$ such that $\ell^m$ kills
the logarithmic class group. As in numerical situations it is rather
easy to compute such an exponent $m$ (when the classical invariants
of the number field are known), this give rise to a more suitable
description of $\wi\Cl_F$ in order to carry on  numerical computations:

%Proposition 2
\begin{Prop}\label{propanni}
Assume the integer $m$ to be large enough such that
the logarithmic class group $\wi\Cl_F$ is annihilated by $\ell^m$. Thus
introduce the group:\smallskip

\centerline{$\wi\Id_F^{(\ell^m)} = \{\a \in \Id_F | \deg_F\d_F(\a) \in 
\ell^m \deg_F \Dl_F \} = \wi\Id_F \ \Id_F^{\ell^m}$}\smallskip

\noindent Thus, denoting $\wi\Pr_F^{(\ell^m)} = \wi\Pr_F \ \wi\Id_F^{\ell^m}$, 
one has: $\wi\Cl_F \simeq \wi\Id_F^{(\ell^m)} / \wi\Pr_F^{(\ell^m)}$ .
\end{Prop}

%\noindent {\it Proof.} 
\begin{proof}
The hypothesis gives $\wi\Id_F^{\ell^m}
\subset \wi\Pr_F$ and by  a straightforward calculation:
$$\begin{aligned}
\wi\Id_F^{(\ell^m)} / \wi\Pr_F^{(\ell^m)} &= 
\wi\Id_F\Id^{\ell^m} / \wi\Pr_F \Id^{\ell^m} \simeq
\wi\Id_F / (\wi\Id_F \cap \wi\Pr_F \Id^{\ell^m})\\
 &\simeq\wi\Id_F / \wi\Pr_F \wi\Id_F^{\ell^m} = \wi\Id_F / \wi\Pr_F
\simeq \wi\Cl_F.
\end{aligned}$$
\end{proof}

\begin{Rem}\label{lemprec0}
\rm
 A lower bound for $m$ which will be required for a sufficient precision of the $p$-adic calculations will be given after lemma \ref{lemprec1}.
\end{Rem}

%%%%%%%%%%%%%%%%%%%%%%%%%%%%%%%%%%%%%%%%%%%

%%%%%%%%%%%%%%%%%%%%%%%%%%%%%%%%%%%%%%%%%%%

%%%%%%%%%%%%%%%%%%%%%%%%%%%%%%%%%%%%%%
%%%%%%%%%%%%%%%%%%%%%%%%%%%%%%%%%%%%%%

\section{The algorithms}

%%%%%%%%%%%%%%%%%%%%%%%%%%%%%%%%%%%%%%
%%%%%%%%%%%%%%%%%%%%%%%%%%%%%%%%%%%%%%

Throughout this section a finite abelian group $G$ is presented
by a column vector $g\in G^m$, whose
entries form a system of generators for $G$, and by a matrix
of relations $M \in \Z^{n\times m}$ of rank $m$, such that
$v^{T}g=0$ for $v \in \Z^m$ if and only if $v^{T}$ is an integral
linear combination of the rows of $M$. We note that for every
$a\in G$ there is a $v\in \Z^m$ satisfying $a=v^{T}g$.
If $g_1, \ldots,g_m$ is a basis of $G$, $M$ is usually a diagonal matrix.
Algorithms for calculations with finite abelian groups can be
found in \cite{cohen-advanced}.
If $G$ is a multiplicative abelian group, then $v^{T}g$ is an
abbreviation for $g_1^{v_1} \cdots g_m^{v_m}$.

%%%%%%%%%%%%%%%%%%%%%%%%%%%%%%%%%%%%%%%%%%%%%%%%%%%%%%%%%%%%%%%%%%

One of the steps in the computation of the logarithmic class group 
is the computation of the ideal class group of a number field.
Algorithms for this can be found in \cite{cohen,florian,PoZa}.
One tool used in these algorithms are the $\s$-units, which we will also 
use directly in our algorithm.

\begin{Def}[$\s$-units]
\rm
Let $\s$ be an ideal of a number field $F$.  We call the group
\[
\bigl\{\alpha\in F^{\times}\mid v_\p(\alpha)=0 \mbox{ for all } \p\nmid\s\bigr\}
\]
the $\s$-units of $F$.
\end{Def}

For this section let $F$ be a fixed number field.  
We denote the ideal class group of $F$ by $\Cl=\Cl_F$.
We also write $\wi\Cl$ for $\wi\Cl_F$, $\wi\Dl$ for $\wi\Dl_F$, and so on.
%%%%%%%%%%%%%%%%%%%%%%%%%%%%%%%%%%%%%%%%%%%%%%%%%%%%%%%%%%%%%%%%%%%%

\subsection{Computing 
\mathversion{bold}
$\deg_F(\p)$ and $\wi v_\p(\cdot)$
\mathversion{normal}}

We describe how invariants of the logarithmic objects can 
be computed.  Some of the tools presented here also applied directly in
the computation of the lo\-garithmic class group.

\begin{DProp}\label{dprop5}
%[J1]
Let $p$ be a prime number.
Let $F$ be a number field.  Let $\p$ be a prime ideal of $F$ over $p$.
For $a \in \Q_p^{\times}\cong  p^\Z\times \F_p^{\times}\times(1+2p\Z_p)$ 
denote by $\langle a \rangle$ the projection of $a$ to
$(1+2p\Z_p)$. 
Let $F_\p$ be the completion of $F$ with respect to $\p$.
For $\alpha\in F$ define
\[
h_\p(\alpha):=\frac{\Log_p \langle\norm_{F_\p/\Q_p}(\alpha)\rangle      }{[F_\p:\Q_p]\cdot\deg_p{p}}.
\]
The $p$-part of the logarithmic ramification index $\wi e_\p$ is
$[h_\p(F_\p^\times):\Z_p]$.  For all primes $q$ with $q\ne p$ the $q$-part of
$\wi e_\p$ is the $q$-part of the ramification index $e_\p$ of $\p$.
\end{DProp}
For a proof see \cite{J1}.

In section \ref{subseclogram} we have seen that the degree $\deg_F(\p)$ of a place 
$\p$ can be computed as  $\deg_F(\p)=\wi f_\p \deg_\l p$.
From section \ref{subsecideal} we know that $\wi e_\p \wi f_\p=e_\p f_\p$. 
We have
\[
\wi v_\p(x) = -\frac{\Log_\l(\norm_{F_\p/\Q_p}(x))}{\deg_F(\p)}.
\]
Thus we can concentrate on the computation of $\wi e_\p$ for which  
we need the completion $F_\p$ of $F$ at $\p$ and generators of the unit group $F_\p^\times$.

The Round Four Algorithm was  originally 
conceived as an algorithm for computing integral bases of number
fields. It can be applied in three different ways in the computation of  
logarithmic classgroups.  Firstly, it is used for factoring
 ideals over number fields; secondly, it returns generating 
 polynomials of completions of number fields; and thirdly, it can be used for 
 determining integral bases of maximal orders.

Let $\Phi(x)$ be a monic, squarefree polynomial over $\Zp$.
The algorithm for factoring polynomials over local fields as described in 
\cite{pauli} returns
\begin{itemize}
\item a factorization $\Phi(x) = \Phi_1(x)\cdot\dots\cdot\Phi_s(x)$
      of  $\Phi(x)$ into irreducible factors $\Phi_i(x)$ $(1\le i\le s)$ over $\Zp$,
\item the inertia degrees $e_i$ and ramification indexes $f_i$ of the extensions 
      of $\Qp$ given by the $\Phi_i(x)$ $(1\le i\le s)$, and
\item two element certificates $(\Gamma_i(x),\Pi_i(x))$ with $\Gamma_i(x),\Pi_i(x)\in F[x]$
      such that 
      $v_i\bigl(\Pi_i(\alpha_i)\bigr)=1/e_i$ and 
      $[\Fp(\ov{\Gamma_i(\alpha_i)}):\Fp]=f_i$ where 
      $\alpha_i$ is a root of $\Phi_i(x)$ in $F[x]/(\Phi_i(x))$,
      $v_i$ is an extension of the exponential valuation $v_p$ of $\Qp$ to $\Qp[x]/(\Phi_i(x))$
      with $v_i|_{\Qp}=v_p$.
\end{itemize}

The factorization algorithm in \cite{ford-pauli-roblot} returns the 
certificates combined in one polynomial for each irreducible factor.  
The data returned by these algorithms can be applied in several ways.
\begin{itemize}
\item
An integral basis of the extension of $\Qp$ generated by a root $\alpha_i$ of
$\Phi_i(x)$ is given by the elements $\Gamma_i(\alpha_i)^h\Pi_i(\alpha_i)^j$ with
$0\leq h\leq f_i$ and $0\leq j \leq e_i$.  The local integral bases 
can be combined to a global integral basis for the extension of
$\Q$ generated by $\Phi(x)$.
\item
For the computation of $\wi v_\p$ we need to 
compute the norm of an element in the completions of $F$.  The completions of $F$ are 
given by the irreducible factors of the generating polynomial of $F$ over $\Q$.
\end{itemize}
%Ideal factorization is a another application.
\begin{Lem}[Ideal Factorization]
Let $\Phi(x)\in\Z_p[x]$ be irreducible over $\Q$.
Let $\Phi_1(x),\dots,\Phi_s(x)\in\Z_p[x]$ be the irreducible factors of $\Phi(x)$ with
two element certificates $(\Gamma_i(x),\Pi_i(x))$.
Denote by  $e_i$  the ramification indexes 
of the extensions
of $\Q_p$ given by the $\Phi_i(x)$ $(1\le i\le s)$.
The Chinese Remainder Theorem gives polynomials $\Theta_1(x),\dots,\Theta_s(x)\in \Q_p[x]$
with
\begin{eqnarray*}
\Theta_i(x)&\equiv&\Pi_i(x)\bmod \Phi_i(x) \\
\Theta_i(x)&\equiv&1 \bmod\textstyle\prod_{j\ne i}\Phi_j(x).
\end{eqnarray*}
Let $L:=\Q(\alpha)$ where $\alpha$ is a root of $\Phi(x)$ in $\CC$.  Then
\[
(p)=\bigl(p,\Theta_1(\alpha)\bigr)^{e_1}\cdot\dots\cdot\bigl(p,\Theta_s(\alpha)\bigr)^{e_s}
\]
is a factorization of $(p)$ into prime ideals.
\end{Lem}

In order to compute $[h_\p(F_\p^\times ):\Z_p]$ it is sufficient to compute the 
image of a set of generators of $F_\p^\times$.
Algorithms for this task were recently
developed with respect to the computation
of ray class groups of number fields and function fields
\cite{cohen-advanced,hepapo}, also see \cite[chapter 15]{hasse}.

\begin{Prop}
\[
F_\p^\times \cong\pi^\Z\times(\OOp/\p)^\times\times(1+\p)
\]
\end{Prop}

Let $\p$ be the prime ideal over the prime number $p$ in $\OOp$.  
Let $e_\p$ be the ramification index and $f_\p$ the inertia
degree of $\p$.  We define the set of fundamental levels
\[
\Fe := \bigr\{\nu  \mid 0 <\nu < {\mbox{$\frac{pe_\p}{p-1}$}}, 
p \nmid \nu \bigl\}
\]
and let $\varepsilon\in\OOp^\times$ such that $p = -\pi^e\varepsilon$.
Furthermore we define the map\smallskip

\centerline{$h_2 : a + \p \longmapsto a^p -\varepsilon a + \p$.}

\begin{Th}[Basis of $(1+\p)$]
Let $\omega_1, \dots,\omega_f\in\OOp$ be a fixed set of representatives of a
$\Fp$-basis of  $\OOp/\p$.
If $(p-1)$ does not divide $e$ or $h_2$ is an isomorphism, then the elements
\[
1+\omega_i\pi^\nu \mbox{ where } \nu \in \Fe, 1\le i\le f\
\]
are a basis of 
the group of principal units $1+\p$.
\end{Th}

\begin{Th}[Generators of $(1+\p)$]
Assume that $(p-1)\mid e$ and $h_2$ is not an isomorphism.
Choose $\enul$ and $\munul$ such that  $p$ does not divide $\enul$ and
such that $e=p^{\munul-1}(p-1)\enul$.
Let $\omega_1, \ldots,\omega_f\in\OOp$ be a fixed set of representatives of a $\Fp$-basis
of $\OOp/\p$ subject to $\omega_1^{p^\munul} -
\varepsilon\omega_1^{p^{\munul-1}} \equiv 0 \bmod \p$.
Choose $\omega_*\in\OOp$ such that $x^p-\varepsilon x
\equiv \omega_* \!\bmod \p$ has no solution.
Then the group of principal units $1+\p$ is generated by
\[
1+\omega_*\pi^{p^{\munul} \enul} \mbox{ and }
1+\omega_i\pi^\nu \mbox{ where } \nu\in\Fe,\,1\le i\le f.
\]
\end{Th}

%%%%%%%%%%%%%%%%%%%%%%%%%%%%%%%%%%%%%%%%%%%%%%%%%%%%%%%%%%%%%%%%%%

\subsection{\mathversion{bold}
Computing a bound for the exponent of $\wi\Cl$
\mathversion{normal}}

%%%%%%%%%%%%%%%%%%%%%%%%%%%%%%%%%%%%%%%%%%%%%%%%%%%%%%%%%%%%%%%%%%

Let $F$ be a number field and $\l$ a prime number.  Let 
$\wi\Cl=\wi\Cl_F\cong \wi\Id/\wi\Pr$ be the $\l$-group of logarithmic divisor classes.
Let $\p_1,\dots,\p_s$ be the $\l$-adic places of $F$.

We describe an algorithm which returns an upper bound $\l^m$ of the exponent of
$\wi\Cl$ (see proposition \ref{propanni}).  We denote by
\begin{itemize}
\item $\wi\Cl(\l)$ the $\l$ group of logarithmic divisor classes of degree zero:
\[
\wi\Cl(\l):=\bigl\{[\a]\in\wi\Cl \mid 
\a=\textstyle\sum_{i=1}^s a_{i} \p_i \mbox{ with }\deg_F(\a)=0\bigr\}
\]
\item $\Cl'$ the $\l$-group of the $\l$-ideal classes, {\sl i.e.}, the $\l$-part
of $\Cl/\!\langle[\p_1],\dots,[\p_s]\rangle$.
\end{itemize}

\begin{Rem}\label{remprec}
\rm
If $(\l)=\p^e$ where $\p$ is a prime ideal of $\OOk$ then the group $\wi\Cl(\l)$ is trivial.
\end{Rem}

%%%%%%%%%%%%%%%%%%%%%%
\begin{Lem}[\cite{DS}]\label{lemds1}
Let
\[\textstyle
\theta:\wi\Cl\longrightarrow\Cl',\;
\sum_\p m_\p\p\longmapsto\prod_{\p\nmid\l}\p^{(1/\lambda_\p) m_\p}.
\]
The sequence
\[
0
\longrightarrow
\wi\Cl(\ell)
\longrightarrow
\wi\Cl
\stackrel{\theta}{\longrightarrow}
\Cl'
\longrightarrow
\Coker\theta
\longrightarrow
0
\]
is exact.
\end{Lem}

\begin{proof} 
%(see \cite{DS})
Recall that, if $\p\nmid\l$,
$\wi v_\p = \lambda_\p v_\p$.
Denote by $\p_1,\dots,\p_s$ the $\l$-adic places of $F$.
Let 
\[
\wi\a=\sum_{\q}a_\q\q=\div(\alpha)
=\sum_\p \wi v_\p(\alpha)\p =
\sum_{\p\nmid\l} \lambda_\p v_\p(\alpha)\p
+\sum_{i=1}^s\wi v_{\p_i}(\alpha)\p_i
\]
be a principal logarithmic 
divisor.  A representative of the image of $\wi\a$ under $\theta$ in terms of 
ideals is of the form 
\[
\a=\prod_{q\nmid(\l)}\q^{v_\q(\alpha)}
=(\alpha\OOk)\times\prod_{i=1}^s \p_i^{-v_{p_i}(\alpha)}.
\]
This shows that the homomorphism $\theta$ is well defined.  
It follows immediately that $\Ker\theta=\wi\Cl(\l)$.
\end{proof}

%%%%%%%%%%%%%%%%%%%%%%

\begin{Lem}\label{lemprec1}
Set $\l^{m'}=\exp\Cl'$ and $\l^{\wi m}=\exp\wi\Cl(\l)$.
Then
\[
\l^{m'+{\wi m}}\a\equiv0\bmod\wi\Pl \mbox{ for all } \a\in\wi\Dl.
\]
\end{Lem}

\begin{proof}
It follows from the exact sequence in lemma \ref{lemds1} that for all
$\a\in\wi\Dl$ the congruence $\l^{m'}\theta(\a)\equiv 1$ holds in $\Cl'$.
Thus $\l^{m'}\a\in\Ker\theta=\wi\Cl(\l)$ and 
$\l^{m'+\wi m}\a\equiv 0  \bmod\wi\Pl $.
\end{proof}

%%%%%%%%%%%%%%%%%%%%%%%

Lemma \ref{lemprec1} suggests setting the precision for the computation of 
$\wi\Cl$ to $m:=m'+\wi m$ $\l$-adic digits.  
If the ideal class group $\Cl$ is known we can easily
compute $m'$.  In order to find $\wi m$ we compute a matrix of relations for 
$\wi\Cl(\l)$.

%%%%%%%%%%%%%%%%%%%%%%%%

\begin{Lem}[Generators and Relations of $\wi\Cl(\l)$]\label{lemprec2}
Let $\p_1,\dots,\p_s$ be the $\l$-adic places of $F$.  Assume that $s>1$.
Reorder the $\p_i$ such that $v_\ell (\deg(\p_1))=\min_{1\le i \le
s}v_\ell (\deg(\p_i))$.  Let  $\gamma_1,\dots,\gamma_r$ be a basis of the
$\l$-units of $F$. Then the group $\wi\Cl(\l)$ is given by the  generators 
$[\g_i]:=\bigl[\p_i-\frac{\deg(\p_i)}{\deg(\p_1)}\p_1\bigr]$ $(i=2,\dots,s)$
with relations 
$\sum_{i=2}^s \wi v_{\p_i}(\gamma_j)[\g_i]=[0]$.
\end{Lem}

\begin{proof}
We consider a logarithmic divisor $\a=\sum_{i=1}^s a_i\p_i$ of degree zero over $F$ 
that is 
constructed from the $\l$-adic places.  
By the choice of $\p_1$ and as
$\deg(\a)=\deg(\sum_{i=1}^s a_i\p_i)=\sum_{i=1}^s a_i\deg(\p_i)=0$
the coefficient $a_1$ is given by the other $s-1$ coefficients.
Thus the $[\g_i]$ generate $\wi\Cl(\l)$.

The relations between the classes of $\wi\Cl(\l)$ are of the form 
$\sum_{i=2}^s b_i [\g_i]=[0]$. 
That is there exists $\beta\in\R_F$ such that
\[
\sum_{i=2}^s b_i \g_i=\sum_{j=1}^s a_j \p_j = \div(\beta),
\]
with 
$v_{\q}(\beta)=0$ for all $\q\nmid(\l)$.
Thus $\beta$ is an element of the group of $\l$-units $\{\alpha\in \R_F \mid v_\q(\alpha)=0\}$
of $\R_F$.  Hence we obtain the relations given above.
\end{proof}

%\begin{Rem}
A version of 
this lemma 
for the case that $F$ is Galois
can be found in \cite{DS}.
%Dieser Beweis ist im galoischen Fall in [DS] schon gemacht !
%\end{Rem}

%%%%%%%%%%%%%%%%%%%%%%%%
\Algo{Precision}
     {algprec}
     {a number field $F$ and  a prime number $\ell$, 
      the $\l$-adic places $\p_1,\dots,\p_s$\\ of $F$, and
      a basis $\gamma_1,\dots,\gamma_r$ of the $\l$-units of $F$} 
     {an upper bound for the exponent of $\wi\Cl$}
{
\begin{itemize}
\itemb Set $\l^{m'}\gets\exp\Cl'$, set $m\gets\max\{m',4\}$.
%\itemb Let $\p_1,\dots,\p_s$ be the $\l$-adic places of $F$.
\itemb If $s=1$ then return $\l^{m'}$.\hfill[Remark \ref{remprec}]
  \itemb Repeat
     \begin{itemize}
     \itemb Set $m\gets m+2$
%     \itemb Let $\gamma_1,\dots,\gamma_r$ be a basis of the $\l$-units of $F$.
     \itemb Set \hfill[Lemma \ref{lemprec2}]
           \[ 
           A \gets \left(\begin{array}{ccc}
                {\wi v_{\p_2}(\gamma_1)} & \dots & {\wi v_{\p_s}(\gamma_1)} \\
                \vdots                   & \ddots& \vdots \\
                {\wi v_{\p_2}(\gamma_r)} & \dots & {\wi v_{\p_s}(\gamma_r)} \\
                \end{array}\right)\bmod\l^m.
           \]
     \itemb Let $H$ be the Hermite normal form of $A$ modulo $\l^m$. 
     \end{itemize}
  \itemb Until $\rank(H)=s-1$.
  \itemb Let $S=(S_{i,j})_{i,j}$ be the Smith normal form of $A$ modulo $\l^m$.
  \itemb Set $\wi m\gets\max_{1\le i\le s-1}\bigl(v_\l(S_{i,i})\bigr)$, 
         return $\l^{m'+\wi m}$.
\end{itemize}
}

\begin{Rem}
\rm
Algorithm \ref{algprec} does not terminate in general if Gross' conjecture
is false.
\end{Rem}
%%%%%%%%%%%%%%%%%%%%%%%%%%

%%%%%%%%%%%%%%%%%%%%%%%%%%%%%%%%%%%%%%

\subsection{
\mathversion{bold}
Computing $\wi\Cl$
\mathversion{normal}
}

%%%%%%%%%%%%%%%%%%%%%%%%%%%%%%%%%%%%%%

We use the ideal theoretic description from section
\ref{subsecideal} for the computation
of $\wi\Cl\cong\wi\Id/\wi\Pr$. 
In the previous section we have seen how we can compute a bound for the 
exponent of $\wi\Cl$.  It is clear that this bound also gives a lower bound for the
precision in our calculations.

%%%%%%%%%%%%%%%%%%%%%%%%%%

\begin{Th}[Generators of $\wi\Cl$]\label{theoclgen}
Let $\a_1,\dots,\a_t$ be a basis of the ideal classgroup 
of $F$ with $\gcd(\a_i,\l)=1$ for all $1\le i \le t$.
Denote by $\p_1,\dots,\p_s$ the $\l$-adic places of $F$.
Let $\alpha_1,\dots,\alpha_s$ be elements of $\R _F$ with
$\wi v_{\p_i}(\alpha_j) = \delta_{i,j}$  $(i,j=1,\dots,s)$ 
and $\gcd((\alpha_i),\l)=1$ for all $1\le i \le s$.
Set $\a_{t+i}:=(\alpha_i)$ for $1\le i \le s$.
%Assume that the $\p_i$ are ordered such that
%$v_\l(\deg(\p_1))=\min_{1\le i \le s}v_\l \deg(p_i)$.
%Set $\beta_i := \alpha_i\cdot\alpha_1^{\deg(\p_i)/\deg(\p_1)}$.
For an ideal $\a$ of $F$ denote by $\ov\a$ the projection of $\a$ from
$\bigoplus_\p \p^{\Z_\l}$ to $\bigoplus_{\p\nmid(\l)} \p^{\Z_\l}$.
We distinguish two cases:
\begin{itemize}
\item[\bf I.] If $\deg_\l(\a_i)=0$ for all  $1\le i \le t+s$ then
set $\b_i:= \a_i$.  The group $\wi\Cl_F$ is generated by 
$\ov\b_1,\dots,\ov\b_{t+s}$.
\item[\bf II.] Otherwise let $1\le j\le t+s$ such that
$v_\l(\deg_\l(\a_j)) =\min_{1\le i\le t+s}v_\l(\deg_\l(\a_i))$.
Set $\b_i := \a_i/\a_j^d$ with
$d \equiv {\frac{\deg_\l(\a_i)}{\deg_\l(\a_j)}} \mod \l^m$ 
where $\l^m > \exp(\wi\Cl)$.
The group $\wi\Cl_F$ is generated by 
$\ov\b_1,\dots,\ov\b_{j-1},\ov\b_{j+1},\dots,\ov\b_{t+s}$.
\end{itemize}
\end{Th}

\begin{proof}
Let $\a\in\wi\Id$.
There exist $\gamma\in\R_F$ and $a_1,\dots,a_t\in\Z_\l$ such that
$\a=\prod_{i=1}^t \a_i^{a_i}\cdot(\gamma)$.
Set $g_i:=\wi v_{\p_i}(\gamma)$ for $1\le i \le s$.
Now
\[
\textstyle\a=\prod_{i=1}^s \a_i^{a_i}
\cdot\bigl((\gamma)\cdot\prod_{j=1}^s(\alpha_i)^{-g_i}\bigr)
\cdot\bigl(\prod_{j=1}^s(\alpha_i)^{g_i}\bigr).
\]
By the definition of $\Id$ (Definition and Proposition \ref{defpropcl}) we have
\[
\textstyle\a=\ov \a=\prod_{i=1}^t \ov \a_i^{a_i}
\cdot\bigl(\ov{(\gamma)\cdot\prod_{j=1}^s(\alpha_j)^{-g_j}}\bigr)
\cdot\bigl(\prod_{j=1}^s\ov{(\alpha_j)^{g_j}}\bigr)
\]
As $\wi v_{\p_i}\bigl((\gamma)\prod_{j=1}^s(\alpha_j)^{-g_j}\bigr)=0$
for $i=1,\dots,s$
we obtain
\[
\textstyle\a\equiv 
\prod_{i=1}^t \a_i^{a_i}
\cdot\bigl(\prod_{j=1}^s{\ov{(\alpha_j)}^{g_j}}\bigr)
\bmod\wi\Pr.
\]
Thus all elements of $\wi\Cl$ can be represented by
$\ov\a_1,\dots,\ov\a_t,\ov\a_{t+1}=\ov{(\alpha_1)},\dots,\ov\a_{t+s}=\ov{(\alpha_s)}$.
For the two cases we obtain:
\begin{itemize}
\item[\bf I.]
It follows immediately that 
$\ov\b_1,\dots,\ov\b_{t+s}$
are generators of $\wi\Cl$.
\item[\bf II.]
If we have 
$
\ov\a \equiv \ov\a_1^{a_1}\cdot\dots\cdot\ov\a_{t+s}^{a_{t+s}}\bmod\wi\Pr
$
for an ideal $\a\in\wi\Id$ 
then  
$
0=\deg(\ov\a)=\sum_{i=1}^{t+s}a_i\deg_\l(\ov\a_i),
$ 
thus
$
-a_j=\sum_{i\ne j}^s a_i\deg_\l(\ov\a_i)/\deg_\l(\ov\a_j).
$
Hence $\ov\b_1,\dots,\ov\b_{j-1},\ov\b_{j+1},\dots,\ov\b_{t+s}$ are generators of $\wi\Cl$.
\end{itemize}
\end{proof}

%%%%%%%%%%%%%%%%%%%%%%%

We continue to use the notation from theorem \ref{theoclgen}.
Set $\Cl':=\Cl/\langle\p_1,\dots,\p_s\rangle$.

\begin{Rem}\rm
The definition of $\Cl'$ in this section and the previous section, 
where we considered the $\l$-part of $\Cl/\langle\p_1,\dots,\p_s\rangle$, differ.
The definition we chose here makes the description of the algorithm easier.
In the algorithm we make sure that only the $\l$-part of the group appears in the
result by computing the $\l$-adic Hermite normal form of the relation matrix.
\end{Rem}

The relations between the generators $\ov\a_1,\dots,\ov\a_t$ of the group $\Cl'$
are of the form $\prod_{i=1}^t\ov\a_i^{a_i}=\ov{(\alpha)}$ with $\alpha\in \R_F$.
There exist integers $c_1,\dots,c_n$ such that 
$\ov{(\alpha)}\equiv \prod_{i=1}^s \ov{(\alpha_i)}^{c_i} \bmod \wi\Pr$.
This yields the relation
$\prod_{i=1}^t\ov\a_i^{a_i}\equiv \prod_{i=1}^s \ov{(\alpha_i)}^{c_i} \bmod \wi\Pr$
in $\wi\Cl$.  
We can derive all relations involving the generators $\ov\a_i+\wi\Pr$ 
from their relations as generators of the group $\Cl'$ in this way.

The other relations between the generators of $\wi\Cl$ are obtained as follows.
A relation between the generators $\ov\alpha_i$ is of the form
$\prod_{i=1}^s\ov{(\alpha_i)}^{v_i}\equiv(1) \bmod\Pr$ or equivalently
$\prod_{i=1}^s(\alpha_i)^{v_i}\cdot\prod_{i=1}^s\p_i^{w_i}=(\alpha)$ 
for some $\alpha\in \R_F$.
The last equality is fulfilled if and only if $\prod_{i=1}^s\p_i^{w_i}$ 
is principal,
{\sl i.e.}, if $\prod_{i=1}^s\p_i^{w_i}$ is an $(\l)$-unit.  
Assume that $\prod_{i=1}^s\p_i^{w_i}=(\gamma)$ for some 
$\gamma\in \R_F$.  
As $\wi v_{\p_j}(\alpha)=0$ for all $\ov{(\alpha)}\in\wi\Pr$ and $\p_j\mid (\l)$ 
the equation 
$\wi v_{\p_j}\bigl(\prod_{i=1}^s\alpha_i^{v_i}\cdot\gamma\bigr)=0$ 
must hold.
By the definition of the $\beta_i$ we obtain $v_i=-\wi v_{\p_i}(\gamma)$
for $1\le i\le s$.

%%%%%%%%%%%%%%%%%%%%%%%%%%

\begin{Cor}[Relations of \mathversion{bold}$\wi\Cl$\mathversion{normal}]\label{corclrel}
Let $\bigl((\ov\a_1,\dots,\ov\a_t), (a_{i,j})_{i,j\in\{1,\dots,t\}}\bigr)$
be a basis and a 
relation matrix of  $\Cl':=\Cl\langle\p_1,\dots,\p_s\rangle$.
Let $\a_{t+1}=(\alpha_1),\dots,\a_{t+s}=(\alpha_s)$ be as above.
For each $1\le k \le t$ we find $c_{k,2},\dots,c_{k,s}$ such that
$\prod_{i=1}^t \ov\b_i^{a_{k,i}}=\prod_{i=2}^s\ov{(\alpha_i)}^{c_{k,i}}$.
Let $\gamma_1,\dots,\gamma_r$ be a basis of the $(\l)$-units of $\R_F$. 
Set $v_{i,j}:=\wi v_{\p_j}(\gamma_i)$ $(1\le i\le r, 2\le j\le s)$.
Set 
\[
M := 
\left(
\begin{array}{cccccc}
b_{1,1} & \dots & b_{1,t} & -c_{1,2} & \dots & -c_{1,s} \\
\vdots  & \ddots& \vdots  & \vdots  & \ddots& \vdots \\
b_{t,1} & \dots & b_{t,t} & -c_{t,2} & \dots & -c_{t,s} \\
  0     & \dots &   0     & v_{1,2} & \dots & v_{1,s} \\
\vdots  & \ddots& \vdots  & \vdots  & \ddots& \vdots \\
  0     & \dots &   0     & v_{r,2} & \dots & v_{r,s} 
\end{array}
\right).
\]
For the two cases we obtain:
\begin{itemize}
\item[\bf I.]$((\ov\b_1,\ov\b_{t+s}),M)$ are generators and relations of $\wi\Cl$.
\item[\bf II.] Let $j$ be chosen as in Theorem \ref{theoclgen}.  
Denote by $N$ the matrix obtained by removing the $j$-th column from $M$.
Then
$((\ov\b_1,\dots,\ov\b_{j-1},\ov\b_{j+1},\dots,\ov\b_{t+s}),N)$ are generators and
relations of $\wi\Cl$.
\end{itemize}
%The group $\wi\Cl$ is given by the generators 
%$\bigl(\ov\b_1,\dots,\ov\b_t,\,\ov{(\beta_2)},\dots,\ov{(\beta_s)}\bigr)$
%and the relation matrix
\end{Cor}
%%%%%%%%%%%%%%%%%%%%%%%%%%

Now we only need to find the elements $\alpha_1,\dots,\alpha_s$ with
$\wi v_{\p_i}(\alpha_j)=\delta_{i,j}$.
Let $\eta_{i,1},\dots,\eta_{i,r_i}$ be a 
system of generators of $\OOp^\times$ for $1\le i \le s$.  Let 
\[
M := 
\left(
\begin{array}{ccc}
\wi v_{\p_1}(\eta_{1,1}) & \dots  &  v_{\p_s}(\eta_{1,1}) \\
\vdots                   & \ddots &  \vdots \\
\wi v_{\p_1}(\eta_{1,r_1}) & \dots  &  v_{\p_s}(\eta_{1,r_1}) \\
\vdots                   & \vdots &  \vdots \\
\wi v_{\p_1}(\eta_{s,1}) & \dots  &  v_{\p_s}(\eta_{s,1}) \\
\vdots                   & \ddots &  \vdots \\
\wi v_{\p_1}(\eta_{s,r_s}) & \dots  &  v_{\p_s}(\eta_{s,r_s})
\end{array}
\right).
\]
Let $S=LMR$ be the $\l$-adic Smith normal form of $M$ with transformation matrices
$L$ and $R$.
Application of the left transformation matrix $L$ to the generators 
$\eta_{1,1},\dots,\eta_{s,r_s}$ yields elements 
$\alpha_1,\dots,\alpha_s$ with the desired properties.

%%%%%%%%%%%%%%%%%%%%%%%%%%

\Algo{Logarithmic Classgroup}
     {alglogcl}
     {a number field $F$ and a prime number $\ell$}
     {generators $g$ and and a relation matrix $H$ for $\wi{Cl}_F$}
{
\begin{itemize}
\item Determine a bound $\l^m$ for the exponent of $\wi{Cl}_F$ and use it as the
      precision for the rest of the algorithm.
      \hfill[Algorithm \ref{algprec}]
\item Compute generators ${\a}_1,\dots,{\a}_t$ of 
      $\Cl'=\Cl/\langle \p_1,\dots,\p_s\rangle$,
      where $\p_1,\dots,\p_s$ are the ideals of $F$ over $\l$.
\item Determine $\a_{t+1}=(\alpha_1),\dots,\a_{t+s}=(\alpha_s)$ with
$\wi v_{\p_i}(\alpha_j)=\delta_{i,j}$.
\item Compute generators $g:=(\ov\b_1,\dots,\ov\b_{t+s})^T$ with $\deg(\b_i)=0$ 
\hfill[Theorem \ref{theoclgen}]
from ${\a}_1,\dots,{\a}_{t+s}$. 
\item Compute a relation matrix  $M$ between the generators $g$. 
      \hfill[Corollary \ref{corclrel}]
\item In case {\bf II.} remove the $j$-th column from $M$ and the $j$-th generator from $g$.
\item Compute the $\l$-adic Hermite normal form $H$ of $M$.
\item Return $(g,H)$.
\end{itemize}
}

%%%%%%%%%%%%%%%%%%%%%%%%%%%%%%%%%%%%%%%%%%%%%%%%%%%%%%%%%%%%%%%%%%%%%%%%%%%%%%%%%%%%%%%%%%%

%%%%%%%%%%%%%%%%%%%%%%%%%%%%%%%%%%%%%%%%%%%%%%%%%%%%%%%%%%%%%%%%%%%%%%%%%%%%%%%%%%%%%%%%%
%%%%%%%%%%%%%%%%%%%%%%%%%%%%%%%%%%%%%%%%%%%%%%%%%%%%%%%%%%%%%%%%%%%%%%%%%%%%%%%%%%%%%%%%%

%%%%%%%%%%%%%%%%%%%%%%%%%%%%%%%%%%%%%%%%%%%%%%%%%%%%%%%%%%%%%%%%%%%%%%%%%%%%%%%%%%%%%%%%%
%%%%%%%%%%%%%%%%%%%%%%%%%%%%%%%%%%%%%%%%%%%%%%%%%%%%%%%%%%%%%%%%%%%%%%%%%%%%%%%%%%%%%%%%%

%%%%%%%%%%%%%%%%%%%%%%%%%%%%%%%%%%%%%%%%%%%%%%%%%%%%%%%%%%%%%%%%%%%%%%%%%%%%%%%%%%%%%%%%%
%%%%%%%%%%%%%%%%%%%%%%%%%%%%%%%%%%%%%%%%%%%%%%%%%%%%%%%%%%%%%%%%%%%%%%%%%%%%%%%%%%%%%%%%%

\section{Examples}

All methods presented here have been implemented in the computer algebra
system Magma \cite{magma}.

We recomputed the logarithmic class groups from 
\cite[section 6]{DS} with our new algorithm.
Our results differ in one example.
For the field $F=\Q(i,\sqrt{1173})$ and $\l=2$ 
we obtain $\wi\Cl_F\cong C_2\times C_2\times C_2$ instead of
$\wi\Cl_F\cong C_2\times C_2\times C_2\times C_2$.
As $F$ contains the 4th roots of unity, the 2-rank of the wild kernel
of $F$ is $3$.

The table contains 
examples of logarithmic $\l$-class groups $\wi\Cl$ of selected number fields $F$
together with their class groups $\Cl$, Galois groups $\Gal$, and the factorization 
of the ideals $(\l)$.
$\chi_\alpha(x)$ denotes the minimal polynomial of $\alpha$ and $i$ 
denotes a root of $x^2+1$.
The class groups are presented as a list of the orders of their cyclic factors,
$\Cl'=\Cl/\langle\p_1,\dots,\p_s\rangle$, and $\l^m$ is the bound for the exponent
of $\wi\Cl$ as returned by algorithm \ref{algprec}.

The logarithmic $2$-class group of $\Q(i,\sqrt{78})$ is an example for the fact that
the cokernel of $\theta$ in the exact sequence in lemma \ref{lemds1} 
is not trivial in general.  Indeed one can show \cite{dubois-soriano} 
that for $F=\Q(i,\sqrt{d})$ with
$d\ne2$ and $d$ squarefree 
\[
\Coker(\theta) \cong
\left\{
\begin{array}{ll}
C_2 & \mbox{ if } d\equiv \pm 2 \bmod 16, \\
C_1 & \mbox{ otherwise.}
\end{array}
\right.
\]

%\begin{table}\label{ourtable}
\begin{center}
\begin{tabular}{llc|rllr|l}

\vpX
$F$ & $\Cl$ & $\Gal$ & 
$\l$ & $(\l)$ & $\Cl'$ & $\!\!\l^m$ & $\wi\Cl$ \\
\hline
\hline

\vpX
$\Q(\sqrt{-521951})$ & [1024] & $S(2)$ &
2 & $\p_1\p_2$ & [4] & $8$ & [2,4]\\
\hline

\vpX
$\Q(i,\sqrt{11})$ & [1] & $E(4)$ &
$5$ &$\p_1\cdots\p_4\!\!\!\!\!\!$  &  [1] & 5 &  [5] \\
\hline

\vpX
$\Q(i,\sqrt{78})$ & [2,2] & $E(4)$ &
$2$ & $\p_1^4$&  [2] & 2 &  [1] \\
\hline

\vpX
$\Q(i,\sqrt{455})$ & [2,2,10] & $E(4)$ &
$2$ & $\p_1^2\p_2^2$&  [2,2] & $\!\!\!\!\!\!\!\!512$ &  [2,512] \\
\hline

\vpX
$\Q(i,\sqrt{1173})$ & [2,2,6] & $E(4)$ &
$2$ & $\p_1^2$&  [2,2,2] & 2 &  [2,2,2] \\
\hline

\vpX
$\Q(i,\sqrt{1227})$ & [4,4] & $E(4)$ &
$\!\!\!\!613$ & $\p_1\cdots\p_4\!\!\!\!\!\!$ & [4,4] & $\!\!\!\!\!\!613$ & [613]\\ 
\hline

\vpX
$\Q(\alpha)$ & [14] & $D(4)$ &
$2$ &$\p_1^2\p_2^2$ &[1]& 1& [1] \\
\multicolumn{3}{l|}{$\chi_\alpha(x)=x^4 + 13x^2 - 12x + 52$}&
\vpX
$3$ &$\p_1^2\p_2^2$ &[1] & 3 &[3]\\
&&&
\vpX
$7$ & $\p_1$ & [14] & 7 & [7]\\
\hline

\vpX
$\Q(\sqrt{1234577},\sqrt{-3})\!\!\!\!\!$ & [273] & $E(4)$ &
$2$  & $\p_1\p_2$ & [273] & 4 & [4,4] \\
\vpX
&&&
$3$  & $\p_1^2$    & [273] & 3 & [3] \\
\vpX
&&&
$13$ & $\p_1\p_2$ & [273] & $\!\!\!\!\!\!169$ & [13,13] \\
\hline
\vpX
$\Q(\zeta_3,\sqrt{303})$ & [14] & $E(4)$ &
$2$ & $\p_1^2$ & [14]  & 2 & [2] \\
\vpX
&&&
$3$ & $\p_1^2\p_2^2$ & [1] & 9 & [9] \\
\vpX
&&&
$7$ & $\p_1\cdots\p_4\!\!\!\!\!\!$ & [1]  & 1 & [1] \\
\hline

%\vpiX
%$\Q(\alpha)$ & [12] & $S(5)$ &
%$2$ &$\p_1\p_2\p_3^3$ &[1]& 1& [1] \\
%\multicolumn{3}{l|}{$\chi_\alpha(x)\!=x^5\!\!+\!2x^4\!\! +\!12*x^3\!\! +14x^2\!\! -12x\! -\!16\!\!$}&
%\vpX
%$3$ &$\p_1\!\!\!\!\!\!$ &[12] & 3 &[3]\\
%\hline

\vpX
$\Q(\beta)$ & [2,6,6] & $S(5)$ &
$2$ &$\p_1\p_2$ &[2,2,6]& 2& [2,2,2] \\
\multicolumn{3}{l|}{$\chi_\beta(x)\!=\!x^5\!\!+\!2x^4\!\!+\!18x^3\!\!+\!34x^2\!\!+\!17x\!+\!3^{10}\!\!$}&
\vpX
$3$ &$\p_1\cdots\p_4\!\!\!\!\!\!$ &[6] & 3 &[3]\\
\hline
\vpX

$\Q(\zeta_5,\sqrt{5029})$ & [15,150] & [2,4] &
$2$ & $\p_1\p_2$ & [3,150]  & 4 & [2,2] \\
\vpX
&&&
$3$ & $\p_1\p_2$ & [15,150]\!\!\!\! & 3 & [3,3] \\
\vpX
&&&
$5$ & $\p_1\p_2$ & [3,150]  & $\!\!\!\!25$ & [5,25] \\
\hline
\vpX

$\Q(i,\sqrt{11},\sqrt{-499})$ &[3,105]& $E(8)$ &
$5$ &$\p_1\cdots\p_8\!\!\!\!\!\!$ & [3] & $\!\!\!\!25$ & [5,5,25]\!\!\!\!\!\!\!\! \\
\hline
\vpX

$\Q(i,\sqrt{11},\gamma)$ &[2,2,2,6]\!\!\!\!& $S(3)\times$ &
$2$ & $\p_1^2$ & [2,2,2,6] & 2 & [2,2,2,2]\!\!\!\!\!\!\!\!\\
\multicolumn{2}{l}{$\chi_\gamma(x)=x^3 + 3x^2 + 2x + 125$}& $E(4)$ &
\vpX
$3$ & $\p_1\p_2$ & [2,2,2,6]\!\!\!\! & 9 & [3,3]\\ 
\vpX
&&&
$5$ & $\p_1\cdots\p_{12}\!\!\!\!\!\!$ &[2]& 5 & [5,5]\\
\end{tabular}
\end{center}
%\end{table}

%%%%%%%%%%%%%%%%%%%%%%%%%%%%%%%%%%%%%%%%%%%%%%%%%%%%%%%%%%%%%%%%%%%%%%%%%%%%%%%%%%%%%%%%%
%%%%%%%%%%%%%%%%%%%%%%%%%%%%%%%%%%%%%%%%%%%%%%%%%%%%%%%%%%%%%%%%%%%%%%%%%%%%%%%%%%%%%%%%%

%%%%%%%%%%%%%%%%%%%%%%%%%%%%%%%%%%%%%%%%%%%%%%%%%%%%%%%%%%%%%%%%%%%%%%%%%%%%%%%%%%%%%%%%%
%%%%%%%%%%%%%%%%%%%%%%%%%%%%%%%%%%%%%%%%%%%%%%%%%%%%%%%%%%%%%%%%%%%%%%%%%%%%%%%%%%%%%%%%%

%%%%%%%%%%%%%%% BIBLIOGRAPHIE %%%%%%%%%%%%%%%%%

{\small
\def\refname{References}
\newcommand{\etalchar}[1]{$^{#1}$}

 }

%%%%%%%%%%%%%%%%%%%%%%%%%%%%%%%%%%%%%%%%%%%%%
\end{document}